\documentclass{amsart}
\usepackage{amsfonts,amsmath,amssymb,amscd,amsthm,enumerate,epsf,graphicx}

 \newtheorem{thm}{Theorem}[section]
 \newtheorem{cor}[thm]{Corollary}
 \newtheorem{lem}[thm]{Lemma}
 \newtheorem{prop}[thm]{Proposition}
 \theoremstyle{definition}
 
 \theoremstyle{remark}
 \newtheorem{rem}[thm]{Remark}
 \newtheorem{exam}[thm]{Example}

 \numberwithin{equation}{section}

 \renewcommand{\H}{\mbox{H}}
 \newcommand{\F}{\mathbb{F}}
 \newcommand{\Q}{\mathbb{Q}}
 \newcommand{\Z}{\mathbb{Z}}

 \newcommand{\Ext}{\rm{Ext}}
 \newcommand{\Tor}{\rm{Tor}}
 
 \newcommand{\Random}{\mbox{Random}}
 \newcommand{\Sum}{\rm{Sum}}
 \newcommand{\reg}{\rm{reg}}
 \newcommand{\ch}{\rm{char}}
 \newcommand{\GL}{\rm{GL}}
 \newcommand{\HS}{\rm{HS}}
 
 \newcommand{\ini}{\rm{in}}
 \newcommand{\pd}{\rm{pd}}
 \newcommand{\lcm}{\rm{lcm}}
 \newcommand{\bideg}{\rm{bideg}}
 \newcommand{\Ideal}{\mbox{Ideal}}
 \newcommand{\fm}{\mathfrak{m}}

 \newcommand{\mbn}{\mathbb{N}}
 \newcommand{\talpha}{\mathbf{\underline{t}}^{\alpha}}
 \newcommand{\xbeta}{\mathbf{\underline{x}}^{\beta}}
 
 \newcommand{\bt}{\beta}
 \newcommand{\ux}{\mathrm{\underline{x}}}
 \newcommand{\ut}{\mathrm{\underline{t}}}

\begin{document}

\title[Linear resolution of powers of ideals]
 {On linear resolution of powers of an ideal}

\author{ Keivan Borna}


\address{Faculty of Mathematical Sciences and Computer, Tarbiat Moallem University, Tehran,
Iran}

\email{borna@ipm.ir}

\thanks{The author is grateful to Dipartimento di Matematica, Universit\'{a} di
Genova, Italia}




\begin{abstract}
In this paper we give a generalization of a result of Herzog,
Hibi, and Zheng providing an upper bound for regularity of powers
of an ideal. As the main result of the paper, we give a simple
criterion in terms of Rees algebra of a given ideal to show that
high enough powers of this ideal have linear resolution. We apply
the criterion to two important ideals $J,J_{1}$ for which we show
that $J^{k},$ and $J_{1}^{k}$ have linear resolution if and only
if $k\neq 2.$ The procedures we include in this work is encoded in
computer algebra package CoCoA \cite{Cocoa}.
\end{abstract}

\maketitle

\section{Introduction}
Let $S = K[x_{1}, \cdots , x_{r}]$ and let
\begin{align*}
\F : \cdots \rightarrow F_{i} \rightarrow F_{i-1}\rightarrow
\cdots
\end{align*}
be a graded complex of free $S$-modules, with $F_{i} = \sum_{j}
S(-a_{i,j})$. The Castelnuovo-Mumford regularity, or simply
regularity, of $\F$ is the supremum of the numbers $a_{i,j} - i.$
The regularity of a finitely generated graded S-module M is the
regularity of a minimal graded free resolution of M. We will write
$\reg(M)$ for this number. The regularity of an ideal is an
important measure of how complicated the ideal is. The above
definition of regularity shows how the regularity of a module
governs the degrees appearing in a minimal resolution. As Eisenbud
mentions in \cite{E02} Mumford defined the regularity of a
coherent sheaf on projective space in order to generalize a
classic argument of Castelnuovo. Mumford's definition \cite{Mu66}
is given in terms of sheaf cohomology. The definition for modules,
which extends that for sheaves, and the equivalence with the
condition on the resolution used above definition, come from
Eisenbud and Goto \cite{EG84}. Alternate formulations in terms of
$\Tor$, $\Ext$ and local cohomology are given in the following.
Let $I$ be a graded ideal, $\fm = (x_{1},\cdots , x_{r})$ the
maximal ideal of $S$, and $n = \dim(S/I).$ Let
\[
a_{i}(S/I) = \max\{t; \H^{i}_{\fm}(S/I)_{t} \neq 0\}, 0 \leq i
\leq n,
\]
where $\H^{i}_{\fm}(S/I)$ is the $i$th local cohomology module
with the support in $\fm$ (with the convention $\max\,\emptyset =
-\infty$). Then the regularity is the number
\[
\reg(S/I) = \max\{a_{i}(S/I) + i ; 0 \leq i \leq n\}.
\]
Note that $\reg(I ) = \reg(S/I ) + 1.$ We may also compute
$\reg(I)$ in terms of $\Tor$ by the formula
\[
\reg(I) = \max_{k} \{t_{k}(I) - k\},
\]
where $t_{p}(I) := \max\{\mbox{degree of the minimal}\, p-\mbox{th
syzygies of}\, I\}.$ Simply this definition may be rewritten as
\begin{align*}
\begin{split}
\reg(I) =& \max_{i,j} \{j - i:\Tor_{i}(I, k)_{j}\neq 0\},\\
       =& \max_{i,j} \{j-i;\beta_{i,j}(I)\neq 0\}.
\end{split}
\end{align*}
Anyway, from local duality one see that the two ways of expressing
the regularity are also connected termwise by the inequality
$t_{k}(I)-k \geq a_{i}(S/I)+n-k.$ Regularity is a kind of
universal bound for important invariants of graded algebras, such
as the maximum degree of the syzygies and the maximum
non-vanishing degree of the local cohomology modules. One has
often tried to find upper bounds for the Castelnuovo-Mumford
regularity in terms of simpler invariants which reflect the
complexity of a graded algebra like dimension and multiplicity.
Clearly $t_{0}(I^{k})\leq k\,t_{0}(I)$ and one may expect to have
the same inequality for regularity, that is, $\reg(I^{k})\leq
k\,\reg(I)$. Unfortunately this is not true in general. 
However, in \cite{CHT99} Cutkosky, Herzog, and Trung and in
\cite{K00} Kodiyalam studied the asymptotic behavior of the
Castelnuovo-Mumford regularity and independently showed that the
regularity of $I^{k}$ is a linear function for large $k$, i.e.,
\begin{align}\label{e:PREG}
\reg(I^{k}) = a(I)k +b(I),\quad \forall k\geq c(I).
\end{align}
Now assume that $I$ is an equigenerated ideal, that is, generated
by forms of the same degree $d$. Then one has $a(I) = d$ and
hence, $\reg(I^{k+1})-\reg(I^{k}) = d$ for all $k \geq c(I)$.
Hence we have
\begin{align}\label{e:LinearPOW}
\reg(I^{k})=(k-c(I))d+\reg(I^{c(I)}),\quad \forall k\geq c(I).
\end{align}
One says that the regularity of the powers of $I$ jumps at place
$k$ if $\reg(I^{k})-\reg(I^{k-1}) > d.$ In \cite{C03} the author
gives several examples of ideals generated in degree $d$ ($d = 2,
3$), with linear resolution (i.e., $\reg(I)=d$), and such that the
regularity of the powers of $I$ jumps at place $2$, i.e., such
that $\reg(I^{2})
> 2d.$ As it is indicated in \cite{C03}, the first example of such
an ideal was given by Terai. Throughout this paper we use $J$ for
this ideal. Geometrically speaking, this is an example of Reisner
which corresponds to the (simplicial complex of a) triangulation
of the real projective plane $\mathbb{P}^{2}$; see Fig.~1 and
\cite{HB} for more details. Let $R:=K[x_{1},\cdots,x_{6}]$ one has
\begin{align}\label{e:Terai}
\begin{split}
J=(&x_{1}x_{2}x_{3},x_{1}x_{2}x_{4},x_{1}x_{3}x_{5},x_{1}x_{4}x_{6},x_{1}x_{5}x_{6},x_{2}x_{3}x_{6},
x_{2}x_{4}x_{5},x_{2}x_{5}x_{6},x_{3}x_{4}x_{5},\\&x_{3}x_{4}x_{6}).
\end{split}
\end{align}
It is known that $J$ is a square-free monomial ideal whose Betti
numbers, regularity and projective dimension depend on the
characteristic of the base field. Indeed whenever $\ch(K)\neq 2,$
$R/J$ is Cohen-Macaulay (and otherwise not), moreover one has
$\reg(J)=3$ and $\reg(J^{2})=7\,(\text{which is of course}>
2\times 3).$ If $\ch(K)=2,$ then $J$ itself has no linear
resolution. So the following natural question arises:


\subsection*{Question A}
 How it goes on for the regularity of powers of $J$?
\\

\begin{figure}
\centerline{
\includegraphics[width=1.8in, height=1.5in]{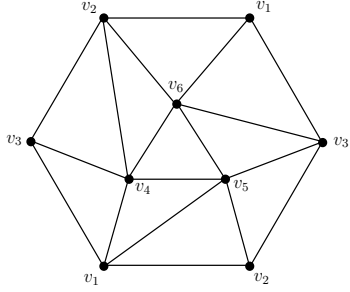}}
\caption{The ideal of triangulation of the real projective plane
$\mathbb{P}^{2}$.}
\end{figure}

By the help of~(\ref{e:PREG}) we are able to write
$\reg(J^{k})=3k+b(J),\,\forall k\geq c(J)$. But what are $b(J)$
and $c(J)$? In this paper we give an answer to this question and
prove that $J^{k}$ has linear resolution (in $\ch(K)=0$) $\forall
k\neq2$, that is, $b(J)=0$ and $c(J)=3$. That is
\begin{align*}
\reg(J^{k})=3k,\quad \forall k\neq 2.
\end{align*}

To answer Question A we develop a general strategy and to this end
we need to follow the literature a little bit. In \cite{R01}
R\"{o}mer proved that
\begin{align}\label{e:Romer}
\reg(I^{n})\leq nd+\reg_{x}(R(I)),
\end{align}
where $R(I)$ is the Rees ring of $I,$ which is naturally bigraded,
and $\reg_{x}$ refers to the x-regularity of $R(I)$, that is,
\[
\reg_{x}(R(I)) = \max\{b - i:\Tor_{i}(R(I), K)_{(b,d)}= 0\},
\]
as defined by Aramova, Crona and De Negri \cite{ACD00}. In Section
~2 we study Rees rings and their bigraded structure in more
details. It follows from ~(\ref{e:Romer}) that if $\reg_{x}(R(I))
= 0$, then each power of $I$ admits a linear resolution. Based on
R\"{o}mer's formula, in ~\cite[Theorem~1.1 and Corollary~1.2]{HHZ}
Herzog, Hibi and Zheng showed the following:

\begin{thm}\label{e:HHZ}
Let $I\subseteq K[x_{1},\cdots,x_{n}]:=S$ be an equigenerated
graded ideal. Let $m$ be the number of generators of $I$ and let
$T:=S[t_{1},\cdots,t_{m}],$ and let $R(I)=T/P$ be the Rees algebra
associated to $I$. If for some term order $<$ on $T,$ $P$ has a
Gr\"{o}bner basis $G$ whose elements are at most linear in the
variables $x_{1},\cdots,x_{n},$ that is $\deg_{x}(f)\leq1$ for all
$f\in G,$ then each power of $I$ has a linear resolution.
\end{thm}

Throughout this paper we simply write $S=K[\ux]$ and $T=S[\ut].$
One can easily see that for $J,$ ~(\ref{e:Terai}), one has at
least ~3 elements in $\ini(P)$ with $deg_{x}>1$, no matter if we
take initial ideal w.r.t. term ordering $\ux>\ut$ or $\ut>\ux$ in
either Lex or DegRevLex order as it is reported in Table ~1. Note
that for example if one starts in DegRevLex order and $\ux>\ut$
then there is ~4 elements in $\ini(P)$ which have ~x-degree $>1$
~($=2\,$ actually) and among them ~2 term has ~t-degree ~1 and ~2
term is in ~t-degree ~2.
\begin{table}
\begin{center}
\begin{tabular}{|c|c|c|}
  \hline
  \, & $\ux>\ut$ & $\ut>\ux$ \\
  \hline
  DegRevLex & ~(1,2):2,(2,2):2 & ~(1,2):2,(2,2):1 \\
  \hline
  Lex & ~(1,2):2,(2,2):1 & ~(1,2):2,(2,2):1 \\
  \hline
\end{tabular}
\end{center}
\caption{Count of elements of $\ini(P)$ with $deg_{x}>1$ for the
ideal of ~(\ref{e:Terai}).}
\end{table}

The main motivation for our work is to generalize Herzog, Hibi and
Zheng's techniques in order to apply them to a wider class.
Furthermore, we will indicate the least exponent $k_{0}$ for which
$I^{k}$ has linear resolution for all $k\geq k_{0}$. Indeed our
generalization works for all ideals which admit the following
condition:

\begin{thm}\label{e:THM}
Let $Q\subseteq S= K[x_{1}, \cdots , x_{r}]$ be a graded ideal
which is generated by m polynomials all of the same degree $d$,
and let $I=\ini(g(P))$ for some linear bi-transformation $g\in
\GL_{r}(K)\times \GL_{m}(K)$. Write $I=G+B$ where $G$ is generated
by elements of $\deg_{x}\leq1$ and $B$ is generated by elements of
$\deg_{x}>1$. If $I_{(k,j)}=G_{(k,j)}$ for all $k\geq k_{0}$ and
for all $j\in\Z$, then $Q^{k}$ has linear resolution for all
$k\geq k_{0}$. In other words, $\reg(Q^{k})=kd$ for all $k\geq
k_{0}$.
\end{thm}

Another motivation for our paper is an example that Conca
considered in \cite{C03}.
\begin{exam}
Let $J_{1}$ be the ideal of ~3-minors of a $4\times4$ symmetric
matrix of linear forms in ~6 variables, that is, ~3-minors of
\[
\left[%
\begin{array}{cccc}
  0 & x_{1} & x_{2} & x_{3} \\
  x_{1} & 0 & x_{4} & x_{5} \\
  x_{2} & x_{4} & 0 & x_{6}\\
  x_{3} & x_{5} & x_{6} & 0 \\
\end{array}%
\right].
\]
As an ideal of $S=\Q[x_{1},\cdots,x_{6}]$ one has:
\begin{align}\label{e:ExDet}
\begin{split}
J_{1}:&=(2x_{1}x_{2}x_{4}, 2x_{1}x_{3}x_{5},2x_{2}x_{3}x_{6},
2x_{4}x_{5}x_{6},x_{1}x_{3}x_{4} + x_{1}x_{2}x_{5} -
x_{1}^{2}x_{6},x_{3}x_{4}x_{6} +\\& x_{2}x_{5}x_{6} -
x_{1}x_{6}^{2},-x_{2}x_{3}x_{4} + x_{2}^{2}x_{5} -
x_{1}x_{2}x_{6},-x_{3}^{2}x_{4} + x_{2}x_{3}x_{5} +
x_{1}x_{3}x_{6},-x_{3}x_{4}^{2} +\\& x_{2}x_{4}x_{5} +
x_{1}x_{4}x_{6},-x_{3}x_{4}x_{5} + x_{2}x_{5}^{2} -
x_{1}x_{5}x_{6}).
\end{split}
\end{align}
\end{exam}
As Conca mentioned in his paper \cite[Remark~3.6]{C03} and as we
will show in this paper, the ideals $J,J_{1}$ are very closely
related. For instance, we prove that
\begin{align*}
\reg(J_{1}^{k})=3k,\quad \forall k\neq 2.
\end{align*}
Similar to the ideal of ~(\ref{e:Terai}), one can easily check
that $\ini(P_{1})$, where $P_{1}$ is the associated ideal to Rees
ring of $J_{1}$, has at least ~9 elements with $deg_{x}>1$, no
matter if we take initial ideal w.r.t. term ordering $\ux>\ut$ or
$\ut>\ux$ in Lex or DegRevLex order; see Table ~2 for more
details. 
\begin{table}
\begin{center}
\begin{tabular}{|c|c|c|}
  \hline
  \, & $\ux>\ut$ & $\ut>\ux$ \\
  \hline
  DegRevLex & ~(1,2):6,(2,2):5,(1,3):1,(4,2):1 & ~(1,2):6,(2,2):3,(1,3):1 \\
  \hline
  Lex & ~(1,2):6,(2,2):3 & ~(1,2):6,(2,2):5 \\
  \hline
\end{tabular}
\end{center}
\caption{Count of elements of $\ini(P_{1})$ with $deg_{x}>1$ for
$J_{1}$,(\ref{e:ExDet}).}
\end{table}

We also show that $J$ and $J_{1}$ and their powers have the same
Hilbert series ($\HS$ for short) correspondingly:
\[
\HS(S/J^{k})=\HS(S/J_{1}^{k}),\quad \forall k.
\]
Indeed we have computed the multigraded Hilbert series of the
corresponding ideals to the Rees algebra of $J$ and $J_{1}$ and
observed that they are the same. As a result we conclude that all
of the powers of $J$ and $J_{1}$ have the same graded Betti
numbers as well:
\[
\bt_{i,j}(J^{k})=\bt_{i,j}(J_{1}^{k}),\quad \forall i,j,\forall k.
\]

\section{Main results}
Let $K$ be a field, $I = (f_{1}, \cdots, f_{m})$ be a graded ideal
of $S = K[x_{1}, \cdots , x_{r}]$ generated in a single degree,
say $d$. The Rees algebra of $I$ is known to be
\[
R(I ) =\bigoplus_{j\geq0} I^{j}t^{j}=S[f_{1}t, \cdots , f_{m}t]
\subseteq S[t].
\]
Let $T = S[t_{1}, \cdots, t_{m}]$. Then there is a natural
surjective homomorphism of bigraded $K$-algebras $\varphi:
T\longrightarrow R(I )$ with $\varphi(x_{i}) = x_{i}$ for $i =
1,\cdots, r$ and $\varphi(y_{j}) = f_{j}t$ for $j = 1,\cdots, m.$
So one can write $R(I) = T/P$. In this paper we consider $T$, and
so $R(I)$, as a standard bigraded polynomial ring with
$\deg(x_{i}) = (0, 1)$ and $\deg(t_{j}) = (1, 0)$. Indeed if we
start with the natural bigraded structure $\deg(x_{i}) = (0, 1)$
and $\deg(f_{j}t) = (d, 1)$ then $R(I)_{(k,vd)}=(I^{k})_{vd}$, but
the standard bidegree normalizes the bigrading in the following
sense:
\begin{align}\label{e:ReesBideg}
R(I)_{(k,j)}=(I^{k})_{kd+j}
\end{align}

For each $k\in \Z$ we define a functor $F_{k}$ from the category
of bigraded $T$-modules to the category of graded $S$-modules with
bigraded maps of degree zero. Let $M$ be a bigraded $T$-module,
define
\[
F_{k}(M)=\bigoplus_{j\in\Z}M_{(k,j)},
\]
obviously $F_{k}$ is an exact functor and associates to each free
$K[{\ux},{\ut}]$-module a free $K[{\ux}]$-module. Sometimes we
simply write $M_{(k,\star)}$ instead of $F_{k}(M)$. Using
~(\ref{e:ReesBideg}) we get
\begin{align}\label{e:fncRees}
[T/P]_{(k,\star)}=R(I)_{(k,\star)}=\bigoplus_{j\in
\Z}R(I)_{(k,j)}=\bigoplus_{j\in \Z}(I^{k})_{kd+j}=I^{k}(kd),
\end{align}
which provides the link between $I$ and its Rees ring $R(I)$. In
the sequel we need to know what is $F_{k}(T(-a,-b))$. For the
convenience of reader we provide a proof.

\begin{rem}
For each integer $k$ we have
\begin{align}\label{e:fncgfree}
T(-a,-b)_{(k,\star)} = \left\lbrace
  \begin{array}{c l}
    0 & \text{if $k<a$},\\
    S(-b)^{N} & \text{otherwise}.
  \end{array}
\right. \end{align} Where $N:=\#\{{\ut}^{\alpha}:
|\alpha|=k-a\}=\binom{m-1+k-a}{m-1}.$
\end{rem}

\begin{proof}
\begin{align}\label{e:NUMERICAL}
\begin{split}
T(-a,-b)_{(k,\star)}&=
\bigoplus_{j\in\Z}T(-a,-b)_{(k,j)}=\bigoplus_{j\in\Z}T_{(k-a,j-b)}\\
&=\bigoplus_{j\in\Z}<\talpha\xbeta:\,|\alpha|=k-a,|\beta|=j-b>,
\end{split}
\end{align}
where the last equality is as vector spaces. From
~(\ref{e:NUMERICAL}) the proof is immediate when $k<a.$
Considering as an $S=K[{\ux}]$-module the last module in
~(\ref{e:NUMERICAL}) is free. Since $|\beta|=j-b$ could be any
integer where $j$ changes over $\Z$, a shift by $-b$ is required
for the representation of the graded free module
$T(-a,-b)_{(k,\star)}$ and finally the proposed $N$ will take care
of the required copies.
\end{proof}

Note that in the spacial case $a=b=0,$ we have
\begin{align}\label{e:fncfree}
T_{(k,\star)}=S^{\binom{m-1+k}{m-1}}.
\end{align}

As we mentioned in Introduction, Theorem ~\ref{e:HHZ} is subject
to condition that $\ini(P)=(u_{1},\cdots,u_{m})$ and
$\deg_{x}(u_{i})\leq 1$. So the natural way to generalize it is to
change the upper bound for $\mbox{x-degree}$ of $u_{i}$ with some
number $t.$ As one may expect, we end up with $\reg(I^{n})\leq
nd+(t-1)\,\pd(T/\ini(P))$. The proof is mainly as that of Theorem
~\ref{e:HHZ} but for the convenience of reader we bring it here.
\begin{prop}\label{e:hhz1}
Let $I\subseteq S$ be an equigenerated graded ideal and let
$R(I)=T/P$. If $\ini(P)=(u_{1},\cdots,u_{m})$ and
$\deg_{x}(u_{i})\leq t$, then $\reg(I^{n})\leq
nd+(t-1)\,\pd(T/\ini(P))$.
\end{prop}
\begin{proof}
Let $C_{\bullet}$ be the Taylor resolution of $\ini(P)$. The
module $C_{i}$ has the basis $e_{\sigma}$ with $\sigma = {j_{1} <
j_{2} < \cdots < j_{i}}\subseteq [m]$. Each basis element
$e_{\sigma}$ has the multidegree $(a_{\sigma}, b_{\sigma})$ where
$x^{a_{\sigma}}. y^{b_{\sigma}} = \lcm\{u_{j_{1}}, \cdots ,
u_{j_{m}}\}$. It follows that $\deg_{x}(e_{\sigma} )\leq ti$ for
all $e_{\sigma}\in C_{i}$. Since the shifts of $C_{\bullet}$ bound
the shifts of a minimal multigraded resolution of \ini(P ), we
conclude that
\begin{align*}
\reg_{x}(T /P)\leq \reg_{x}(T/\ini(P))&=  \max_{i,j}\{a_{ij}-i\}\\
&\leq  ti-i=(t-1)i\\ &\leq (t-1)\,\pd(T/\ini(P)).
\end{align*}
Now ~(\ref{e:Romer}) completes the proof.
\end{proof}
One can see that now Theorem ~\ref{e:HHZ} is the special case of
Proposition ~\ref{e:hhz1} with $t=1.$ However, this approach seems
to be less effective. Our approach to generalize Theorem
~\ref{e:HHZ} is to change
$P$ with an isomorphic image $g(P)$ 
so that $\ini(g(P))_{(k,\star)}$ only consists of 
terms with ~x-degree$\leq$~1, for some $k$. To this end, we need a
simple fact.

Let $<$ be any term order on $S=K[\ux]$ and let $V\subseteq S$ be
a $K$-vector space. Then with respect to the monomial order on $S$
obtained by restricting $<$, by definition $V$ is homogenous if
for any element $f$ of $V$, $f=\displaystyle\sum_{i=0}^{n} f_{i}$,
where $f_{i}$ is an element of $S$ of degree $i$, we have
$f_{i}\in V,\,\forall i=0,\cdots,n.$ That is to say $
V=\displaystyle\bigoplus_{i=0}^{\infty} V_{i}\,,\,V_{i}=V\bigcap
S_{i}.$ It yields that
$\ini(V)=\displaystyle\bigoplus_{i=0}^{\infty} \ini(V_{i})$ and
so, $\ini(V)_{i}=\ini(V_{i}).$ Generalizing this idea to bigraded
(or multigraded) situation is also well understood. Let $F$ be a
free $S$-module with a fixed basis and $M$ a bigraded subvector
space of it. Then
\[
\ini(M)_{(i,j)}=\ini(M_{(i,j)}),
\]
and so
\begin{align}\label{e:iniFnc1}
\ini(M)_{(k,\star)}:=\displaystyle\bigoplus_{j\in\Z}\ini(M)_{(k,j)}=
\displaystyle\bigoplus_{j\in\Z}\ini(M_{(k,j)})=\ini(M_{(k,\star)}).
\end{align}
See \cite{E01} chapter ~15.2 for more details. Furthermore since
$\bt_{ij}^{S}(F/M)\leq\bt_{ij}^{S}(F/\ini(M)),$ it is easy to
conclude with
\begin{align}\label{e:regini}
\reg(F/M)\leq\reg(F/\ini(M)).
\end{align}

\begin{lem}\label{e:regIniFnc}
Let $P$ be the associated ideal of Rees ring $R(I)$ and let
$T=R/P$. Then
$\reg([T/P]_{(k,\star)})\leq\reg([T/\ini(P)]_{(k,\star)}).$

\begin{proof}Since $P$ is a naturally bigraded ideal of T, and since
easily $T_{(k,\star)}$ is a free $S$-module (see
~(\ref{e:fncfree})), ~(\ref{e:iniFnc1}) implies that
$\ini(P)_{(k,\star)}=\ini(P_{(k,\star)}).$
Applying ~(\ref{e:regini}) for $F:=T_{(k,\star)}$ and $M:=P$ we
obtain $\reg(T_{(k,\star)}/P_{(k,\star)})
\leq\reg(T_{(k,\star)}/\ini(P_{(k,\star)})).$ Finally putting all
together we get the required inequality.
\begin{align*}\label{e:iniFnc3}
\begin{split}
\reg([T/P]_{(k,\star)})=\reg(T_{(k,\star)}/P_{(k,\star)})
&\leq\reg(T_{(k,\star)}/\ini(P_{(k,\star)}))\\
&=\reg(T_{(k,\star)}/\ini(P)_{(k,\star)})\\
&=\reg([T/\ini(P)]_{(k,\star)}).
\end{split}
\end{align*}
\end{proof}
\end{lem}

In the following the proof of Theorem ~\ref{e:THM} is given.
\begin{proof}
First of all notice that, since $g:K[\ux,\ut]\longrightarrow
K[\ux,\ut]$ is an invertible bi-homogenous transformation, we have
the following bi-homogenous isomorphism:
\[
\frac{K[\ux,\ut]}{P}\simeq  \frac{K[\ux,\ut]}{g(P)},
\]
and so we can simply take $g=id$ in the rest of proof. 
Write down the so-called Taylor resolution of $T/G$: 
\begin{align}\label{e:TJ1}
\begin{array}{ccccccccccc}
                                                                          \! & \! & F_{2,0}  & \! & \! & \! & \! & \! & \! & \! & \! \\
                                                                    \! & \! & \bigoplus & \! & F_{1,0} & \! & \! & \! & \! & \! & \! \\
  \cdots&\longrightarrow & F_{2,1} & \longrightarrow & \bigoplus & \longrightarrow & T & \longrightarrow & T/G & \longrightarrow & 0, \\
                                                                    \! & \! & \bigoplus & \! & F_{1,1} & \! & \! & \! & \! & \! & \! \\
                                                                           \! & \! & F_{2,2} & \! & \! & \! & \! & \! & \! & \! & \! \\
\end{array}
\end{align}
where $F_{i,j}=\bigoplus_{a\in\Z} T(-a,-j)^{\bt_{i,(a,j)}(T/G)}.$
Note that $\bt_{i,(a,j)}(T/G),$ is an integer number which depends
on $i$, $a$, and $j$. Since $(k,\star)$ is an exact functor, the
following complex of $K[\ux]$-modules is exact:
\begin{align}\label{e:TJ2}
\begin{array}{ccccccccccc}
                                                                          \! & \! & (F_{2,0})_{(k,\star)}  & \! & \! & \! & \! & \! & \! & \! & \! \\
                                                                    \! & \! & \bigoplus & \! & (F_{1,0})_{(k,\star)} & \! & \! & \! & \! & \! & \! \\
  \cdots&\longrightarrow & (F_{2,1})_{(k,\star)} & \longrightarrow & \bigoplus & \longrightarrow & T_{(k,\star)} & \longrightarrow & [T/G]_{(k,\star)} & \longrightarrow & 0. \\
                                                                    \! & \! & \bigoplus & \! & (F_{1,1})_{(k,\star)} & \! & \! & \! & \! & \! & \! \\
                                                                           \! & \! & (F_{2,2})_{(k,\star)} & \! & \! & \! & \! & \! & \! & \! & \! \\
\end{array}
\end{align}
Using formula ~(\ref{e:fncgfree}) we obtain
$T(-a,-b)_{(k,\star)}=S(-b)^{N_{a,k}}$, so for $F_{i,j}$ we get
\begin{align}\label{e:TJ3}
(F_{i,j})_{(k,\star)}=\bigoplus_{a\in\Z}
S(-j)^{N_{a,k}\,\bt_{i,(a,j)}(T/G)}.
\end{align}
It follows that ~(\ref{e:TJ2}) is a (possibly non-minimal) graded
free $K[\ux]$- resolution of $[T/G]_{(k,\star)}$. Since
$\deg_{x}(G)\leq1$, from ~(\ref{e:TJ2}) and ~(\ref{e:TJ3}) we
conclude that
\begin{align}\label{e:regTG}
 \reg([T/G]_{(k,\star)})=0 \quad \text{for all $k$}.
\end{align}
Now we have
\begin{align}\label{e:MENQ}\begin{split}
dk\leq \reg(Q^{k})\leq \reg([T/P]_{(k,\star)})+dk & \leq
\reg([T/\ini(P)]_{(k,\star)})+dk \\ & = \reg([T/G]_{(k,\star)})+dk \quad\text{for all $k\geq k_{0}$}\\
&=0+dk=dk,
\end{split}\end{align}
 where the second (in)equality
in ~(\ref{e:MENQ}) follows from ~(\ref{e:fncRees}), the third
inequality is due to Lemma ~\ref{e:regIniFnc}, and the forth comes
from the easy argument
$[T/\ini(P)]_{(k,\star)}=T_{(k,\star)}/\ini(P)_{(k,\star)}=T_{(k,\star)}/G_{(k,\star)}=[T/G]_{(k,\star)}.$

Finally ~(\ref{e:MENQ}) implies that $\reg(Q^{k})=kd$ for all
$k\geq k_{0}$ as desired.
\end{proof}


\section{Examples and applications}
In this section we provide some applications of Theorem
\ref{e:THM}. But before that we examine our condition on the
decomposition of $\ini(P)$ in a closer view. In the following a
reformulation of our results is provided.

With the assumptions and notation introduced in Theorem
~\ref{e:THM} assume that $B=(m_{1},\cdots,m_{p})$ and
$\bideg(m_{i})=(t_{i},\geq2).$ By $(t_{i},\geq2)$ we mean that the
$\deg_{x}(m_{i})\geq 2$. It is harmless to assume that
$t_{1}\leq\cdots\leq t_{p}.$ If for all $i=1,\cdots,p$ and all
$\alpha\in\mbn^{m}$ with $\mid\alpha\mid=t_{p}+1-t_{i}$ we have
$\ut^{\alpha}m_{i}\subseteq G$ then $I_{(k,\star)}=G_{(k,\star)}$
for all $k>t_{p}+1.$

Using this strategy and as an application for our main result we
give an answer to the Question A proposed in the Introduction. 

\begin{exam}\label{e:Traie}
Let $S=\Q[x_{1},\cdots,x_{6}]$ and let $J$ be the ideal of
(\ref{e:Terai}). Let
$T=\Q[x_{1},\cdots,x_{6},t_{1},\cdots,t_{10}]$ with order
$\ux>\ut$ (and DegRevLex). We also use $J$ for the ideal of $T$
generated by the same generators as of $J$ in $S$. Let $P$ be the
defining ideal of the Rees ring of $J$, so $R(J)=T/P$. One can
check that $P$ has ~15 elements of bidegree ~(1,1), ~10 elements
of bidegree ~(3,0), and ~15 elements of bidegree ~(4,0). Take $G$
and $B$ as in Theorem ~\ref{e:THM}.
We have checked that $|G|=60,\,B=\Ideal(t_{6}x_{4}x_{5},
t_{4}x_{3}x_{5}, t_{4}t_{6}x_{5}^{2}),$ and so
$\max\{\deg_{t}(h)\mid h\in B\}=2.$ But
$(\ut)^{2}(t_{6}x_{4}x_{5})\nsubseteq G,
(\ut)^{2}(t_{4}x_{3}x_{5})\nsubseteq G,
(\ut)(t_{4}t_{6}x_{5})\nsubseteq G.$
So in DegRevLex (also Lex) order and $\ux>\ut,$ we were unable to
admit the conditions of Theorem ~\ref{e:THM}. We have observed
that the same story happens for ordering $\ut>\ux$ either
DegRevLex or Lex. One could try to take $g$ "generic", as in
~(\ref{e:Rand}).
\begin{align}\label{e:Rand}
\begin{split}
&g:=g_{1}\times g_{2},\\
&g_{1}:=x_{i}\longmapsto\Random(\Sum(x_{1},\cdots,x_{6})),\\
&g_{2}:=t_{j}\longmapsto\Random(\Sum(t_{1},\cdots,t_{10})),
\end{split}
\end{align}
for all $i=1,\cdots,6$ and all $j=1,\cdots,10$, where by
$\Random(\Sum(x_{1},\cdots,x_{6}))$ we mean a linear combination
of $x_{1},\cdots,x_{6}$ with random coefficients and the same
interpretation for $t_{1},\cdots,t_{10}.$ But we realized that a
properly chosen sparse random upper triangular $g$ does the job as
well. We continue in DegRevLex order and $\ut>\ux$.

We have implemented some functions (in CoCoA) to look for a
desired upper triangular bi-change of coordinates. For example,
the following $g$ works fine for $J$, indeed there exists many of
such $g$:
\begin{align*}
g:=g_{1}\times g_{2}\in \GL_{6}(\Q)\times \GL_{10}(\Q),
\end{align*}
where $g_{1}:\Q[\ux]\longrightarrow \Q[\ux]$ is given by
\begin{align*}
&x_{4}\longmapsto x_{1}+x_{4},\\
&x_{6}\longmapsto x_{3}+x_{6},
\end{align*}
and sends $x_{i}$ for $i\neq4,6$ to itself and let $g_{2}$ to be
the identity map over $\Q[\ut]$. One can compute that $\mid
G\mid=98,\,B=(t_{7}x_{3}^{2}, t_{4}t_{6}x_{5}^{2}).$ It is easy to
verify that
\begin{align}\label{e:con31}
I_{(k,\star)}=G_{(k,\star)},\,\text{for\;}k> 2\Longleftrightarrow
\left\{\begin{array}{cc}(t_{7}x_{3}^{2})(t_{1},\cdots,t_{10})^{2}\subseteq
G,\\
(t_{4}t_{6}x_{5}^{2})(t_{1},\cdots,t_{10})\subseteq G,
\end{array}\right.
\end{align}
and since in the right side of ~(\ref{e:con31}) both containments
are valid we conclude with $\reg(J^{k})=3k$ for all $k>2.$
\end{exam}

Taking several ideas from Example ~\ref{e:Traie} now we are able
to quickly find an answer to Question A for $J_{1}.$ In the
following we show that $\reg(J_{1}^{k})=3k,$ for all $k>2.$

\begin{exam}\label{e:Ex36}
Let $S=\Q[x_{1},\cdots,x_{6}]$ and let $J_{1}$ be the ideal of
~(\ref{e:ExDet}). Let
$T=\Q[t_{1},\cdots,t_{10},x_{1},\cdots,x_{6}]$ in DegRevLex order,
and let $P_{1}$ be the defining ideal of the Rees ring of $J_{1}$,
so $R(J_{1})=T/P_{1}$. One can observe that $P$ has ~15 elements
of bidegree ~(1,1), ~10 elements of bidegree ~(3,0), and ~12
elements of bidegree ~(4,0). Take $g$ to be the following simple
upper triangular bi-transformation:
\begin{align*}
g:=g_{1}\times g_{2}\in \GL_{6}(\Q)\times \GL_{10}(\Q),
\end{align*}
where $g_{1}:\Q[\ux]\longrightarrow \Q[\ux]$ shall be given by
\begin{align*}
&x_{4}\longmapsto x_{2}+x_{4},\\
&x_{6}\longmapsto x_{1}+x_{6},
\end{align*}
and sending the rest to themselves and take
$g_{2}:\Q[\ut]\longrightarrow \Q[\ut]$ to be
\begin{align*}
t_{8}\longmapsto t_{7}+t_{8},
\end{align*}
and for $i\neq 8,\,t_{i}\longmapsto t_{i}$. Computations by CoCoA
shows that $|G|=144,\,B=(t_{10}x_{2}x_{3}, t_{2}t_{4}x_{5}^{2}).$
Since $I:=\ini(g(P))= G+B$, we have
\begin{align}\label{e:con36}
I_{(k,\star)}=G_{(k,\star)},\,\text{for\;}k> 2\Longleftrightarrow
\left\{\begin{array}{cc}(t_{10}x_{2}x_{3})(t_{1},\cdots,t_{10})^{2}\subseteq
G,\\
(t_{2}t_{4}x_{5}^{2})(t_{1},\cdots,t_{10})\subseteq G,
\end{array}\right.
\end{align}
and since it is easy to check that the right side of
~(\ref{e:con36}) is holding, we obtain that $\reg(J_{1}^{k})=3k$
for all $k>2.$
\end{exam}

We conclude with the following two corollaries which indicate that
ideals $J,$ ~(\ref{e:Terai}), and $J_{1}$, ~(\ref{e:ExDet}), are
very tightly related.
\begin{cor}
When the characteristic of the base field is zero, all the powers
of $J$, and $J_{1}$, but the second power have linear resolution.
\end{cor}

Since the least exponent $k_{0}$ for $J^{k}$, and also for
$J_{1}^{k}$ in order to have linear resolution for all $k>k_{0}$
is ~2, the following question seems to be interesting to discover:

\subsection*{Question B}Does there exist an ideal $Q$ with generators of the same degree $d$ over some polynomial ring
$S=K[x_{1},\cdots,x_{r}]$, for which $\reg(Q^{k})=kd,\,\forall
k\neq 3$ or $\forall k\neq 2,3$?
\\

As we mentioned in Introduction, it is easy to check that $T/P$
and $T/P_{1}$ have the same multigraded Hilbert series, where $P$,
and $P_{1}$ are the defining ideals of Rees rings of $J$ and
$J_{1}$ correspondingly. The immediate result is as follows:
\begin{cor}
$\HS(S/J^{k})=\HS(S/J_{1}^{k})\,\forall k,$ and so
$\bt_{i,j}(J^{k})=\bt_{i,j}(J_{1}^{k})\,\forall i,j,\forall k.$
\end{cor}

\section*{Acknowledgment} The results of this paper were obtained
during the visit of the author to Dipartimento di Matematica,
Universit\'{a} di Genova, Italia. The author would like to express
his deep gratitude to Professor Aldo Conca for the kind invitation
and for his warm hospitality whose guidance and support were
crucial for the successful completion of this project. During the
stay in Genova, the author was supported by a grant within the
frame of the Italian network on ``Commutative, Combinatorial,
Computational Algebra'' PRIN 2006-07, directed by Professor Valla.
It is a pleasure for the author to warmly thank him and also the
kind staff in DIMA as well. Finally, the role of the software
package CoCoA in computations of concrete examples as we worked on
this project is acknowledged.

\bibliographystyle{amsplain}

\end{document}